\documentclass[12pt]{article}
\usepackage{amsmath,latexsym}
\usepackage{amsthm}
\usepackage{enumerate}
\usepackage{epsfig}

\newcommand\eps{{\varepsilon}}

\newcommand{\thankyou}{Supported in part by Japanese GCOE Program G08: ``Fostering Top Leaders in Mathematics - Broadening the Core and Exploring New Ground".}

\newcommand{\address}{Address: Department of Mathematics, University of North Texas, 1155 Union Circle \#311430, Denton, TX 76203-5017, USA; E-mail: kiko@unt.edu}

\newtheorem{theorem}{Theorem}[section]

\newtheorem{lemma}[theorem]{Lemma}
\newtheorem{remark}[theorem]{Remark}

\title{At which points exactly has Lebesgue's singular function the derivative zero ?}

\author{Kiko Kawamura \\University of North Texas \footnote{\thankyou}\ \,\footnote{\address}}
\date{\today}

\begin{document}
\pagestyle{myheadings}

\maketitle
\begin{abstract}

Let $L_a(x)$ be Lebesgue's singular function with a real parameter $a$ ($0<a<1, a \neq 1/2$). As is well known, $L_a(x)$ is strictly increasing and has a derivative equal to zero almost everywhere. However, what sets of $x \in [0,1]$ actually have $L_a^{'}(x)=0$ or $+\infty$? We give a partial characterization of these sets in terms of the binary expansion of $x$. As an application, we consider the differentiability of the composition of Takagi's nowhere differentiable function and the inverse of Lebesgue's singular function.  



\end{abstract}

\section{Introduction}

Imagine flipping an unfair coin with probability $a \in (0,1)$ of heads and probability $1-a$ of tails. Note that $a \neq 1/2$. 
Let the binary expansion of $t \in [0,1]$: $t=\sum_{n=1}^{\infty}\omega_n / 2^n$ be determined by flipping the coin infinitely many times. 
More precisely, $\omega_n=0$ if the $n$-th toss is heads and $\omega_n=1$ if it is tails. 
We define {\it Lebesgue's singular function}  $L_a(x)$ as the distribution function of $t$:
\begin{equation*}
L_a(x):=prob\{t \leq x\}, \qquad 0 \leq x \leq 1.
\end{equation*}

\begin{figure}
  \begin{center}
    \epsfig{file=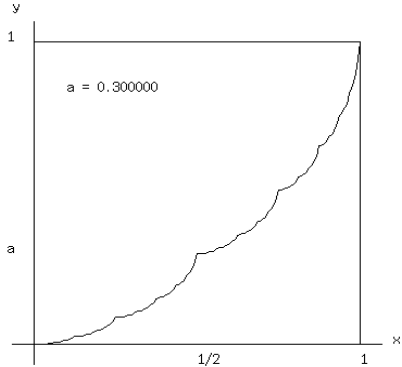,height=2.0in,width=.4\textwidth}
  \end{center}
  \caption{Lebesgue's singular function ($a=0.3$)}
\end{figure}

It is well-known that $L_a(x)$ is strictly increasing, but the derivative is $0$ almost everywhere. 
This distribution function $L_a(x)$ was also defined in different ways and studied by a number of authors: Cesaro (1906), Faber (1910), Lomnicki and Ulam(1934), Salem (1943), De Rham (1957) and others. 
For instance, De Rham~\cite{DeRham-1957} studied $L_a(x)$ as a unique continuous solution of the functional equation
\begin{equation} \label{eq:FE}
        L_{a}(x)=
        \begin{cases}
                a L_{a}(2x), 
                        & \qquad 0 \leq x \leq \tfrac12, \\
                (1-a)L_{a}(2x-1)+a, 
                        & \qquad \tfrac12 \leq x \leq 1, \\
        \end{cases}
\end{equation}
where $0<a<1$, and $a \neq 1/2$. 

From \eqref{eq:FE}, it is clear that the graph of $L_a(x)$ is self-affine. Because of its connection with fractals, several applications have been found in recent years: for instance, in physics~\cite{Takayasu-1984, Tasaki-1993}, real analysis~\cite{Hata-Yamaguti-1984, Kawamura-2002}, digital sum problems~\cite{Kruppel-2009, Okada-Sekiguchi-Shiota-1995} and complex dynamical systems~\cite{Sumi-2009}. There is even a connection with the Collatz conjecture~\cite{Berg-Kruppel-2000}. 

\medskip

Reconsider the differentiability of $L_a(x)$. 
It is known that for any $x \in [0,1]$, $L^{'}_{a}(x)$ is either zero, or $+\infty$, or it does not exist. Then, it is natural to ask at which points $x \in [0,1]$ exactly we have $L^{'}_{a}(x)=0$ or $+\infty$. 

\medskip

In fact, De Rham~\cite{DeRham-1957} gave the following partial answer to this question. Let the binary expansion of $x \in [0,1]$ be $x=\sum_{k=1}^\infty 2^{-k}\eps_k$, where $\eps_k \in \{0,1\}$. For those $x \in [0,1]$ having two binary expansions, we choose the expansion which is eventually all zeros. As an exception, fix $\eps_k=1$ for every $k$ if $x=1$. 

Define 
\begin{equation}\label{eq:In}
I_n:=\sum_{k=1}^{n}\eps_k.
\end{equation} 
Note that $I_n$ is the number of $1$'s occurring in the first $n$ binary digits of $x$. 

Suppose that $I_n/n$ tends to a limit $l$ as $n \to \infty$, and let 
\begin{equation}\label{eq:l_0}
l_0:=\frac{\log 2a}{\log a - \log(1-a)}. 
\end{equation}
Then the derivative $L^{'}_{a}(x)$ exists and is zero, when $(l-l_0)(a-1/2)>0$. 
An English translation of De Rham's original paper is included in Edgar's book \cite{Edgar-1993}.

\medskip

Unfortunately, De Rham's paper did not contain a proof. The main purpose of this note is to give a proof of De Rham's statement and extend his result. The paper is organized as follows. Section 2 states and proves the main results. The key to the proof is to use Lomnicki and Ulam's expression from  1934~\cite{Lomnicki-Ulam-1934}. De Rham might have had a different proof in mind, as he did not mention Lomnicki and Ulam's paper. In Section 3, as an application, we consider a question about the differentiability of the composition of Takagi's nowhere differentiale function and the inverse of Lebesgue's singular function.

\section{The main result}

For convenience, define the right-hand and left-hand derivatives of $L_a(x)$ as follows.
\begin{align*}
L'_{a+}(x):&=\lim_{h \to 0+} \frac{L_a(x+h)-L_a(x)}{h}, \\
L'_{a-}(x):&=\lim_{h \to 0-} \frac{L_a(x+h)-L_a(x)}{h},
\end{align*}
provided the limits exist. 

From the self-affinity of the graph, we have 
\begin{lemma}\label{lem1}
For any $x \in [0,1]$ for which $L'_{a+}(x)$ exists, 
$$L'_{a+}(x)=L'_{(1-a)-}(1-x).$$
\end{lemma} 

\medskip

Define 
\begin{equation} 
D_{1}(x):= \lim_{n \to \infty}\frac{I_n}{n}=\lim_{n \to \infty}\frac{1}{n}\sum_{k=1}^{n}\eps_k,
\end{equation}
provided the limit exists, and put $D_{0}(x):=1-D_{1}(x)$.
In other words, $D_{i}(x)$ is the density of the digit $i$ in the binary expansion of $x$, for $i=0,1$.

\begin{theorem}\label{th:main}
\begin{enumerate}
\item If $x \in[0,1]$ is dyadic, then $L'_{a+}(x) \neq L'_{a-}(x)$.

\item If $x \in [0,1]$ is not dyadic and $0<D_1(x)<1$, then
\begin{equation*}
L'_a(x)=\begin{cases}
0, & \mbox{if }\quad a^{D_0(x)}(1-a)^{D_1(x)}< 1/2,\\
+\infty, & \mbox{if }\quad a^{D_0(x)}(1-a)^{D_1(x)}> 1/2.
\end{cases}
\end{equation*}
\end{enumerate}
\end{theorem}

\begin{remark}
De Rham's statement is equivalent to the following. For a value of $x$ for which $D_1(x)$ exists,  $L'_a(x)=0$ when $a^{D_0(x)}(1-a)^{D_1(x)}< 1/2$.
\end{remark}

\begin{remark}
If $x$ is a binary normal, that is, $D_0(x)=D_1(x)=1/2$, then Theorem~\ref{th:main} gives $L'_a(x)=0$, since $\sqrt{a(1-a)}<1/2$.
\end{remark}

\bigskip
\noindent {\bf Proof of Theorem \ref{th:main}}.
First, suppose $x \in [0,1]$ is a dyadic point, say $x=j/2^N$. 
Let $2^{-(k+1)}\leq h \leq 2^{-k}$ where $k>N$. Since $L_a$ is increasing, this implies that
\begin{equation}\label{eq: inequality-1}
\frac{L_a(x+2^{-(k+1)})-L_a(x)}{2^{-k}} \leq \frac{L_a(x+h)-L_a(x)}{h} \leq \frac{L_a(x+2^{-k})-L_a(x)}{2^{-(k+1)}}.
\end{equation}

The key to the proof is to use the following expression for $L_a(x)$, given by Lomnicki and Ulam~\cite{Lomnicki-Ulam-1934}:
\begin{equation}\label{eq:Ulam}
L_{a}(x)=\frac{a}{1-a}\sum_{n=1}^{\infty}\eps_n a^{n-I_n}(1-a)^{I_n}, 
\end{equation}
where $I_n$ is defined by \eqref{eq:In}. By \eqref{eq:Ulam}, we have 
\begin{equation*}
L_a(x+2^{-k})-L_a(x)=a^{k-I_N}(1-a)^{I_N}.
\end{equation*}
Since $(1-a)/a$ is a positive constant, 
\begin{align*}
\lim_{k \to \infty}\frac{L_a(x+2^{-k})-L_a(x)}{2^{-k}}&=\lim_{k \to \infty}(2a)^k \left(\frac{1-a}{a}\right)^{I_N} \\
&=\begin{cases}
0, & \mbox{if }\quad 0 < a < 1/2,\\
+\infty, & \mbox{if }\quad 1/2 < a < 1.
\end{cases}
\end{align*}
By \eqref{eq: inequality-1}, it follows that 
\begin{equation*}
L'_{a+}(x)=\begin{cases}
0, & \mbox{if }\quad 0 < a < 1/2,\\
+\infty, & \mbox{if }\quad 1/2 < a < 1.
\end{cases}
\end{equation*}

Since $1-x$ is also a dyadic, the left-hand derivative follows from Lemma~\ref{lem1}:
\begin{equation*}
L'_{a-}(x)=L'_{(1-a)+}(1-x)= \begin{cases}
+\infty, & \mbox{if }\quad 0 < a < 1/2,\\
0, & \mbox{if }\quad 1/2 < a < 1.
\end{cases}
\end{equation*}
Therefore, $L'_{a}(x)$ does not exist if $x$ is dyadic.

\bigskip

Next, suppose $x \in [0,1]$ is not dyadic and $0<D_1(x)<1$. 
Let $p_k$ be the address of the $k$-th $``0"$ in the binary expansion of $x$, and $2^{-p_{k+1}} \leq h \leq 2^{-p_k}$. 
Since $L_a$ is increasing, this implies that
\begin{equation}\label{eq: inequality-2}
\frac{L_a(x+2^{-p_{k+1}})-L_a(x)}{2^{-p_{k}}} \leq \frac{L_a(x+h)-L_a(x)}{h} \leq 
\frac{L_a(x+2^{-p_{k}})-L_a(x)}{2^{-p_{k+1}}}.
\end{equation}

Using \eqref{eq:Ulam}, we have 
\begin{equation}\label{eq:la-1}
L_a(x+2^{-p_k})-L_a(x)=a^k(1-a)^{p_k-k}+\left(1-\frac{a}{1-a}\right)\sum_{n=p_k+1}^{\infty}\eps_na^{n-I_n}(1-a)^{I_n}.
\end{equation}
Let $n(l)$ be the address of the $l$-th $``1"$ appearing after position $p_k$ in the binary expansion of $x$. Then we have
\begin{equation*}
\sum_{n=p_k+1}^{\infty}\eps_na^{n-I_n}(1-a)^{I_n}=a^k(1-a)^{p_k-k}\sum_{l=1}^{\infty}a^{n(l)-p_k-l}(1-a)^{l}.
\end{equation*}
Since $n(l)-p_k-l \geq 0$ and $0<a<1$, the series in the right hand side above converges, say to $C(x,k)$. 

For convenience, define 
\begin{equation*}
C_1(x,k):=1+\left(1-\frac{a}{1-a}\right)C(x,k).
\end{equation*}
Then we can write \eqref{eq:la-1} as 
\begin{equation}\label{key-equation}
L_a(x+2^{-p_k})-L_a(x)=a^k(1-a)^{p_k-k}C_1(x,k).
\end{equation}
Since $C(x,k) \leq \sum_{l=1}^{\infty}(1-a)^l$, it follows that  
\begin{equation}\label{eq:c1-bounded}
\min \left\{1, \frac{1-a}{a}\right\} \leq C_1(x,k) \leq \max\left\{1, \frac{1-a}{a}\right\}.
\end{equation}

By \eqref{key-equation}, we have 
\begin{align*}
\frac{L_a(x+2^{-p_k})-L_a(x)}{2^{-p_{k+1}}}
&=\left\{2^{\frac{p_{k+1}}{p_k}} a^{\frac{k}{p_k}}(1-a)^{1-\frac{k}{p_k}}\right\}^{p_k} C_1(x,k), \\
\frac{L_a(x+2^{-p_{k+1}})-L_a(x)}{2^{-p_k}}
&=\left\{2^{\frac{p_k}{p_{k+1}}} a^{\frac{k+1}{p_{k+1}}}(1-a)^{1-\frac{k+1}{p_{k+1}}}\right\}^{p_{k+1}} C_1(x,k).
\end{align*}

Since $k/p_k$ tends to a nonzero limit $D_0(x)$ as $k \to \infty$, we have $p_{k+1}/p_k \to 1$ as $k \to \infty$. 
Therefore, it follows from \eqref{eq: inequality-2} and \eqref{eq:c1-bounded} that 

\begin{equation*}
L'_{a+}(x)=\begin{cases}
0, & \mbox{if }\quad a^{D_0(x)}(1-a)^{D_1(x)}< 1/2,\\
+\infty, & \mbox{if }\quad a^{D_0(x)}(1-a)^{D_1(x)}> 1/2.
\end{cases}
\end{equation*}

\bigskip

Finally, for the left-hand derivative, it follows from Lemma~\ref{lem1} that 
\begin{equation*}
L'_{a-}(x)= L'_{(1-a)+}(1-x)=\begin{cases}
0, & \mbox{if }\quad a^{D_0(x)}(1-a)^{D_1(x)}< 1/2,\\
+\infty, & \mbox{if }\quad a^{D_0(x)}(1-a)^{D_1(x)}> 1/2,
\end{cases}
\end{equation*}
since $D_i(x)=D_j(1-x)$ when $i \neq j$. This concludes the proof. 
$\Box$

\bigskip
\begin{remark}

A careful study of the above proof shows that the existence of the full limit $D_1(x)=\lim_{n \to \infty}(I_n/n)$ is not necesary. The following generalization is straightforward:
\begin{enumerate}
\item Suppose $0<a<1/2$. Then
\begin{equation*} 
L'_{a}(x)=\begin{cases}
0, & \mbox{if }\quad \lim_{n \to \infty}\sup I_n/n < l_0, \\
+\infty, & \mbox{if }\quad \lim_{n \to \infty}\inf I_n/n > l_0.
\end{cases}
\end{equation*}

\item Suppose $1/2<a<1$. Then 
\begin{equation*} 
L'_{a}(x)=\begin{cases}
0, & \mbox{if }\quad \lim_{n \to \infty}\inf I_n/n > l_0, \\
+\infty, & \mbox{if }\quad \lim_{n \to \infty}\sup I_n/n < l_0,
\end{cases}
\end{equation*}
where $l_0$ is defined by \eqref{eq:l_0}.
\end{enumerate}
\end{remark}

\bigskip

Note that Theorem \ref{th:main} left out the boundary case; that is, those numbers $x$ for which $a^{D_0(x)}(1-a)^{D_1(x)}=1/2$; in other words, numbers $x$ which have the following densities: 
\begin{equation*}\label{boundary condition}
D_1(x)=\frac{\log 2a}{\log a -\log(1-a)}, \qquad D_0(x)=\frac{\log 2(1-a)}{\log(1-a)-\log a}.
\end{equation*}

Let us define some additional notation. As a complement of $I_n$, define $O_n$ to be the number of $0$'s occurring in the first $n$ binary digits of $x$: 
\begin{equation*} \label{eq:O_n}
O_n:=\sum_{k=1}^{n}(1-\eps_k). 
\end{equation*} 
Let $q_k$ be the address of the $k$-th $"1"$ in the binary expansion of $x$ as a complement of $p_k$. Observe that 
\begin{align*} 
q_k \leq n \qquad & \mbox{if and only if } \qquad I_n \geq k, \\
p_k \leq n \qquad & \mbox{if and only if } \qquad O_n \geq k. 
\end{align*}

Then, it is easy to prove the following lemma by contradiction. 

\begin{lemma}\label{lem:key}
Let $f(k)=p_k-\frac{k}{D_0(x)}$ and $g(k)=\frac{k}{D_1(x)}-q_k$. 
If $f(k) \to \infty$ as $k \to \infty$, then $g(k) \to \infty$.
\end{lemma}
%
%

\begin{theorem}\label{theorem:boundary}
Suppose $x \in [0,1]$ satisfies $a^{D_0(x)}(1-a)^{D_1(x)}=1/2$. 
Let $f(k)=p_k-\frac{k}{D_0(x)}$ and suppose $f(k+1)/f(k) \to 1$. 
\begin{enumerate}
\item If $f(k) \to \infty$ as $k \to \infty$, then  
\begin{align*}
L'_{a}(x)=
					&\begin{cases}
+\infty, & \mbox{if }\quad 0<a<1/2,\\
0, & \mbox{if }\quad 1/2<a<1.
\end{cases}
\end{align*}
\item If $f(k) \to -\infty$ as $k \to \infty$, then 
\begin{align*}
L'_{a}(x)=
					&\begin{cases}
0, & \mbox{if }\quad 0<a<1/2,\\
+\infty, & \mbox{if }\quad 1/2<a<1.
\end{cases}
\end{align*}
\end{enumerate}
\end{theorem}

\noindent {\bf Proof of Theorem \ref{theorem:boundary}}.

We follow the same argument for non-dyadic points $x \in [0,1]$ as in the proof of Theorem \ref{th:main}.
Let $f(k)=p_k-\frac{k}{D_0(x)}$. Since $k/p_k$ tends to a nonzero limit $D_0(x)$ as $k \to \infty$, $f(k)$ is of smaller order than $k$. 
Then, it follows from \eqref{key-equation} that 
\begin{align*}
\frac{L_a(x+2^{-p_k})-L_a(x)}{2^{-p_k}}&=\left[\left\{2(1-a)\right\}^{\frac{1}{D_0(x)}}a (1-a)^{-1}\right]^k \left\{2(1-a)\right\}^{f(k)}C_1(x,k) \\
&=\left\{2(1-a)\right\}^{f(k)}C_1(x,k), 
\end{align*}
because 
\begin{equation*}
\{2(1-a)\}^{\frac{1}{D_0(x)}}a(1-a)^{-1}=1, \qquad \mbox{when } a^{D_0(x)}(1-a)^{D_1(x)}=1/2.
\end{equation*}

Thus,
\begin{align*}
\frac{L_a(x+2^{-p_k})-L_a(x)}{2^{-p_{k+1}}}&=\left\{2^{\frac{f(k+1)}{f(k)}}(1-a)\right\}^{f(k)}\cdot 2^{\frac{1}{D_0(x)}} C_1(x,k), \\
\frac{L_a(x+2^{-p_{k+1}})-L_a(x)}{2^{-p_{k}}}&=\left\{2^{\frac{f(k)}{f(k+1)}}(1-a)\right\}^{f(k+1)}\cdot 2^{\frac{1}{D_0(x)}} C_1(x,k+1).
\end{align*}

Since $f(k+1)/f(k) \to 1$ as $k \to \infty$, it follows from \eqref{eq: inequality-2} and \eqref{eq:c1-bounded} that if $f(k) \to \infty$, 
\begin{align*}
L'_{a+}(x)
	& =\begin{cases}
+\infty, & \mbox{if }\quad 0<a<1/2,\\
0, & \mbox{if }\quad 1/2<a<1.
\end{cases}
\end{align*}

Similary, if $f(k) \to -\infty$ as $k \to \infty$, then 
\begin{align*}
L'_{a+}(x)
					& =\begin{cases}
0, & \mbox{if }\quad 0<a<1/2,\\
+\infty, & \mbox{if }\quad 1/2<a<1.
\end{cases}
\end{align*}

Next, consider the left-hand derivative. From Lemma \ref{lem1}, we have $L'_{a-}(x)=L'_{(1-a)+}(1-x)$. 
It is clear that $1-x$ also satisfies $a^{D_0(1-x)}(1-a)^{D_1(1-x)}=1/2$, since $D_i(x)=D_j(1-x)$ for $i \neq j$. 
Let $g(k)=\frac{k}{D_1(x)}-q_k$. Since $q_k=p_k(1-x)$, we have $g(k+1)/g(k) \to 1$ if $f(k+1)/f(k) \to 1$. 
It follows from Lemma \ref{lem:key} that if $f(k) \to \infty$ or $-\infty$, then $L'_{a-}(x)=L'_{(1-a)+}(1-x)=L'_{a+}(x)$. This concludes the proof. $\Box$

\section{Application}

%

We apply the main result to the following simple question. 
In classical calculus, the chain rule is used to compute the derivative of the composition of two differentiable functions. However, what can we say, for example, about the differentiability of the composition of a nowhere differentiable function and a singular function? For instance, let $T$ be  Takagi's nowhere differentiable function, which is defined by 
\begin{equation*}
T(x)=\sum_{k=0}^{\infty}\frac{1}{2^{k}}|2^kx-\lfloor 2^kx+\tfrac12 \rfloor|, \qquad 0 \leq x\leq 1. 
\end{equation*}
Is $(T \circ L_a^{-1})$ nowhere differentiable? See Figure $2$. If $a=0.4$, the figure of the graph looks somewhat like Takagi's function; on the other hand, if $a=0.2$, the shape of the graph is more like Lebesgue's singular function. Thus, we can guess that $(T \circ L_a^{-1})$ might not be nowhere differentiable if $a$ is close to $0$. 

\medskip
Although $T$ does not have a finite derivative anywhere, it is known to have an improper infinite derivative at many points. In fact, Allaart and Kawamura~\cite{Allaart-Kawamura-2010} proved that the set of points where $T'(x)=+\infty$ or $-\infty$ has Hausdorff dimension one. Note that the inverse of Lebesgue's singular function is also singular. Hence, if we try to (naively) use the chain rule to compute the derivative of 
$(T \circ L_a^{-1})(x)$, we may run into one of the indeterminate products $+\infty\cdot 0$ or $-\infty\cdot 0$. 

\bigskip

The following theorem gives an answer to this concrete question: $(T \circ L_a^{-1})(x)$ has a finite but vanishing derivative at uncountably many points.
\begin{figure}
\begin{center}
	\epsfig{file=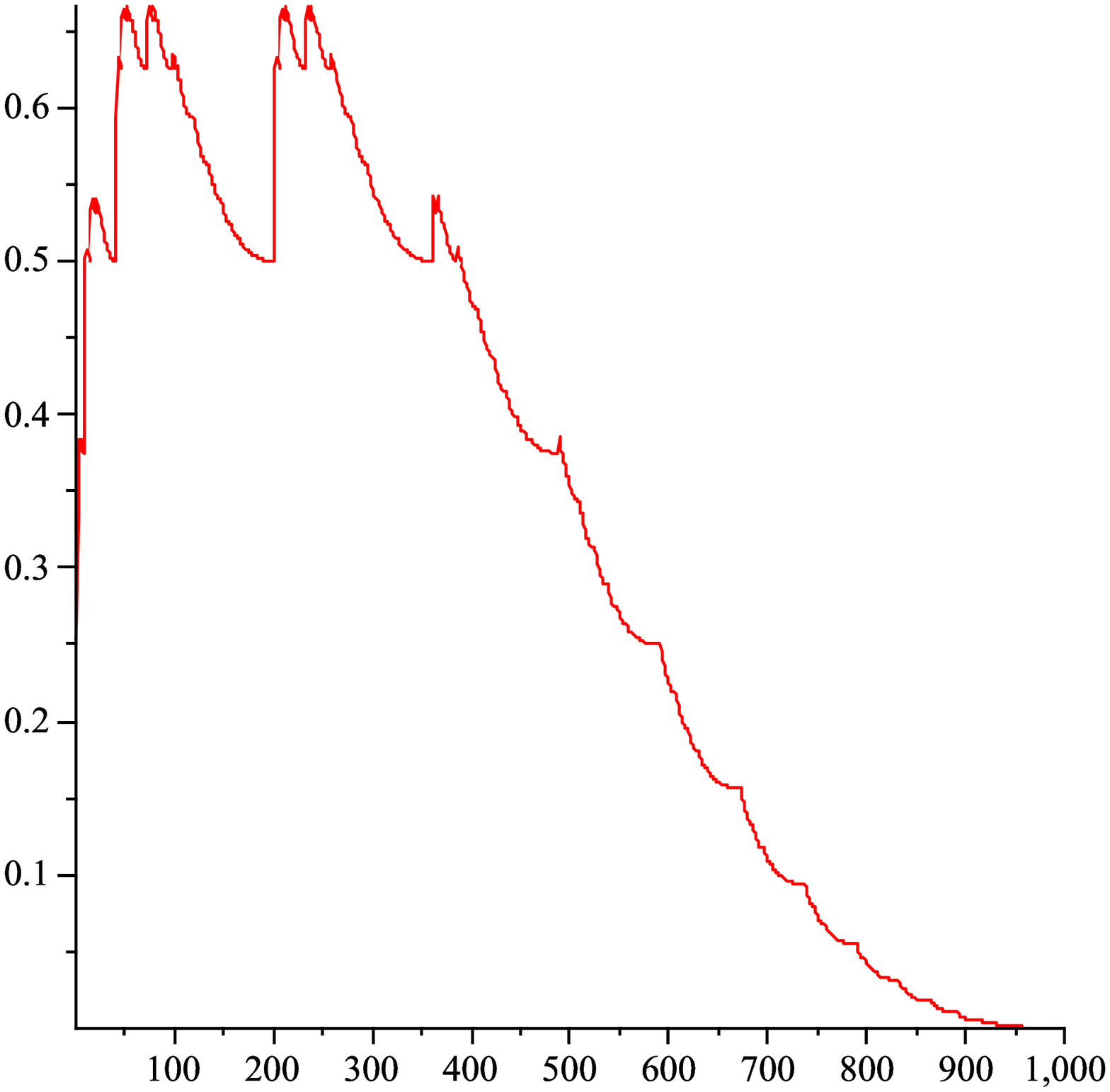,height=1.2in,width=.3\textwidth}
	\epsfig{file=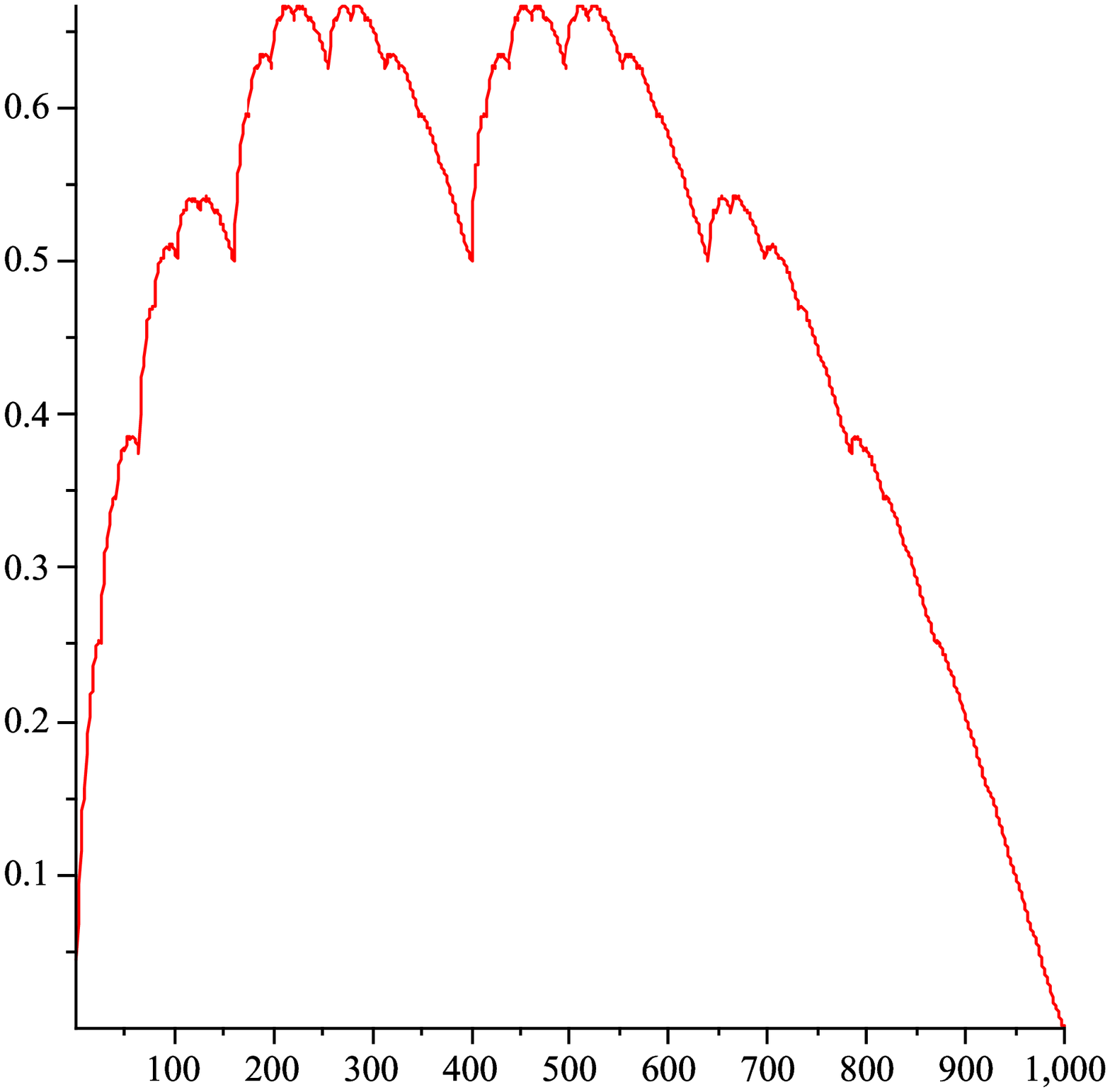,height=1.2in,width=.3\textwidth}
\caption{Graphs of $(T\circ L_a^{-1})(x)$ for $a=0.2$ (left) and $a=0.4$ (right)}
\end{center}
\end{figure}

\begin{theorem}\label{th:main2}
Let $x\in[0,1]$, and put $y=L_a^{-1}(x)$.
If $0<D_1(y)<1$ and $a^{D_0(y)}(1-a)^{D_1(y)}> 1/2$, then 
\begin{equation}\label{th3:composition}
(T\circ L_a^{-1})'(x)=0.
\end{equation}
\end{theorem}

\bigskip
\noindent {\bf Proof}.

\bigskip
Define $\tilde{h}:=L_a^{-1}(x+h)-L_a^{-1}(x)$. Then we can write 
\begin{equation}\label{modulus-composition}
\frac{T(L_a^{-1}(x+h))-T(L_a^{-1}(x))}{h} 
=\frac{T(y+\tilde{h})-T(y)}{\tilde{h}\log_2(1/|\tilde{h}|)} \cdot \frac{\tilde{h}\log_2(1/|\tilde{h}|)}{h}.
\end{equation}

Allaart and Kawamura~\cite{Allaart-Kawamura-2010} proved that if $D_1(x)$ exists and $0<D_1(x)<1$, then
\begin{equation*}
\lim_{h \to 0} \frac{T(x+h)-T(x)}{h\log_2(1/|h|)}=D_{0}(x)-D_{1}(x).
\end{equation*} 

Therefore, we have 
\begin{equation*} 
-1 \leq \lim_{h \to 0}\frac{T(y+\tilde{h})-T(y)}{\tilde{h}\log_2(1/|\tilde{h}|)} \leq 1.
\end{equation*}

A slight modification of the proof of Theorem \ref{th:main} yields 
\begin{equation*}
\lim_{h \to 0}\frac{\tilde{h}\log_2(1/|\tilde{h}|)}{h}=0, \qquad \mbox{if }\quad a^{D_0(y)}(1-a)^{D_1(y)}> 1/2.
\end{equation*}
Substituting these results into \eqref{modulus-composition} gives \eqref{th3:composition}.
$\Box$

\section*{Acknowledgment}

This research was done mainly during my visit to RIMS, Kyoto university. I am grateful to Prof.~H.~Okamoto for his support and warm-hearted hospitality. Also, I would like to thank Prof.~P.~Allaart for his helpful comments and suggestions in preparing this paper. 

Lastly, I wish to dedicate this paper to the memory of Prof.~H.~Shinya, who taught me a deeper understanding of calculus.


\footnotesize
\begin{thebibliography}{10}


\bibitem{Allaart-Kawamura-2010}
{\sc P.~Allaart} and {\sc K.~Kawamura}\@, The improper infinite derivatives of Takagi's nowhere-differentiable function,  {\em J. \ Math. \ Anal. \ Appl.}, {\bf 372}, pp~656-665 \ (2010). 

\bibitem{Berg-Kruppel-2000}
{\sc L.~Berg} and {\sc M.~Kruppel}\@, De Rham's singular function and related functions, 
{\em Z. \ Anal. \ Anwendungen.}, {\bf 19}, no.~1, pp~227-237 \ (2000). 


\bibitem{DeRham-1957}
{\sc G.~de Rham}\@, Sur quelques courbes d\'efinies par des \'equations fonctionnelles, {\em Rend.\ Sem.\ Mat.\ Torino} {\bf 16}, pp.~101-113 \ (1957).

\bibitem{Edgar-1993}
{\sc G.~A.~Edgar}\@, {\it Classics on Fractals}, Addison-Wesley, Reading, MA (1993).


\bibitem{Hata-Yamaguti-1984}
{\sc M.~Hata} and {\sc M.~Yamaguti}\@, Takagi function and its generalization, {\em Japan J.\ Appl.\ Math.},\ {\bf 1}, pp.~183-199 \ (1984).

\bibitem{Kawamura-2002}
{\sc K.~Kawamura}\@, On the classification of self-similar sets determined by two contractions on the plane,
{\em J.\ Math. \ Kyoto Univ.},\ {\bf 42}, no.~2, pp.~255-286 \ (2002).

\bibitem{Kruppel-2009}
{\sc M.~Kruppel}\@, De Rham's singular function, its partial derivatives with respect to the parameter and binary digital sums, {\em Rostock. \ Math. \ Kolloq.}, {\bf 64}, pp~57-74 \ (2009). 


\bibitem{Lomnicki-Ulam-1934}
{\sc Z.~Lomnicki} and {\sc S.~Ulam}\@, Sur la th\'eorie de la mesure dans les espaces combinatoires et son application au calcul des probabilit\'es I. Variables ind\'ependantes. {\em Fund.\ Math.}\ {\bf 23}, pp.~237-278 \ (1934).


\bibitem{Okada-Sekiguchi-Shiota-1995}
{\sc T.~Okada}, {\sc T.~Sekiguchi} and {\sc Y.~Shiota}\@, 
An explicit formula of the exponential sums of digital sums, {\em Japan J. Indust.\ Appl.\ Math.}\ {\bf 12}, pp.~425-438 \ (1995).




\bibitem{Sumi-2009}
{\sc H.~Sumi}, Rational semigroups, random complex dynamics and singular functions on the complex plane. {\em Sugaku}
 {\bf 61}, no. 2, pp.~133-161 (2009).  


\bibitem{Takagi-1903} 
{\sc T.~Takagi}\@, A simple example of the continuous function without derivative, {\em Phys.-Math. Soc. Japan} {\bf 1} (1903), 176-177. {\em The Collected Papers of Teiji Takagi}, S. Kuroda, Ed., Iwanami, pp.~5-6 (1973). 

\bibitem{Takayasu-1984}
{\sc H.~Takayasu}\@, Physical models of fractal functions,
 {\em Japan J.\ Appl.\ Math.},\ {\bf 1}, pp.~201-205 \ (1984).

\bibitem{Tasaki-1993}
{\sc S.~Tasaki}, {\sc I.~Antoniou} and {\sc Z.~Suchanecki}\@, Deterministic diffusion, De Rham equation and fractal  eigenvectors,\ {\em Physics Letter} A 179, pp.~97-102 \ (1993).

\end{thebibliography}
\end{document}